\documentclass{article}

\usepackage{amssymb}

\usepackage{amsfonts,amsmath,amsthm}
\usepackage{graphicx,epstopdf}
\usepackage{graphicx}
\usepackage{float}

\usepackage{listings}
\usepackage{color}

\usepackage{float}
\usepackage{amsmath} 
\usepackage{calc}
\newlength{\depthofsumsign}
\usepackage{calc}

\newcommand{\RR}{\mathbb{R}}

\begin{document}



\newcommand{\Eqref}[1]{(\ref{#1})}


\begin{center}
	\large{ \textbf{ {\Large Approximating basins of attraction for dynamical systems via stable radial bases}}}
\end{center}

\begin{center}
	Roberto Cavoretto$^*$,  Stefano De Marchi$^+$, Alessandra De Rossi$^*$, \\
	Emma Perracchione$^*$, Gabriele Santin$^+$
\end{center}

\begin{center}
	$^*$Department of Mathematics "G. Peano", University of Turin - Italy\\
	$^+$Department of Mathematics, University of Padova - Italy
\end{center}
\vskip 0.5cm

\textbf{Abstract.} 	In applied sciences, such as physics and biology, it is often required to model the evolution of populations via dynamical systems. 
In this paper, we focus on the problem of approximating the basins of attraction of such models in case of multi-stability. 
We propose to reconstruct the domains of attraction via an implicit interpolant using stable radial bases, obtaining the surfaces by partitioning 
the phase space into disjoint regions. An application to a competition model presenting jointly three stable equilibria is considered.

\section{Introduction}
Mathematical modelling is nowadays commonly used in applied sciences, such as in biology, to predict the temporal evolution of populations. 
This is obtained in general via a dynamical system, where a particular solution is completely determined by its initial conditions. 
The possible steady states of the system are determined by its parameters, 
but more steady states can originate from different initial conditions \cite{Arrowsmith90,Murray02}. 

Here, focusing on a particular competition model presenting tristability, we approximate the basins of attraction 
of its stable equilibria. To do that, we propose an efficient and stable kernel-based interpolation method to reconstruct (unknown) multiple surfaces partitioning the three-dimensional space or the phase space into disjoint regions. Since in many situations the standard basis interpolant may suffer from instability due to ill-conditioning of the kernel matrices, we consider a change of basis that makes the related interpolant more stable than the standard one. This computational approach extends our recent research in this field (see \cite{Cavoretto14,Cavoretto15a,Cavoretto15b,Demarchi13,Demarchi15}).

Denoting now the three populations by $x$, $y$ and $z$, we consider the competition model
\begin{align} \label{model3d}
\begin{array}{ll}
\frac{ \displaystyle  dx}{ \displaystyle  dt}=p \big(1- \frac{ \displaystyle  x}{ \displaystyle  u} \big)x-axy-bxz,  & \textrm{} \\
\vspace{.01cm}\\
\frac{ \displaystyle  dy}{ \displaystyle  dt}=q \big(1- \frac{ \displaystyle  y}{ \displaystyle  v} \big)y-cxy-eyz, & \textrm{} \\
\vspace{.01cm}\\
\frac{ \displaystyle  dz}{ \displaystyle  dt}=r \big(1- \frac{ \displaystyle  z}{ \displaystyle  w} \big)z-fxz-gyz,  & \textrm{} 
\end{array}
\end{align}
where all the parameters characterizing the system are positive and defined as follows: 
\begin{itemize}
	\item $p$, $q$ and $r$ are the growth rates of $x$, $y$ and $z$, respectively;
	\item $a$, $b$, $c$, $e$, $f$ and $g$ are the competition rates;
	\item   $u$, $v$ and $w$ are the carrying capacities of the three populations.
\end{itemize}
The model \eqref{model3d} has eight equilibria. The origin $E_0 = (0, 0, 0)$ and the points associated with the survival of 
only one population $E_1 = (u, 0, 0)$, $E_2 = (0, v, 0)$ and $E_3 = (0, 0, w)$. Then we have three equilibrim points with two coexisting populations, i.e., 
$$
E_4 =  \bigg( \frac{\displaystyle uq(av-p)}{\displaystyle cuva-pq},\frac{\displaystyle pv(cu-q)}{\displaystyle cuva-pq},0\bigg), 
$$
$$
E_5 = \bigg( \frac{\displaystyle ur(bw-p)}{\displaystyle fuwb-rp},0, \frac{\displaystyle wp(fu-r)}{\displaystyle fuwb-rp}\bigg), 
 $$
 $$
E_6 = \bigg( 0, \frac{\displaystyle vr(we-q)}{\displaystyle gvwe-qr}, \frac{\displaystyle wq(vg-r)}{\displaystyle gvwe-qr}\bigg).
$$
Finally, we have the coexistence equilibrium $E_7=(x_7,y_7,z_7)$, which can be assessed only via numerical simulations.

For example, taking as parameters $p = 1$, $q = 2$, $r =2$, $a = 5$, $b=4$, $c =3$, $e=7$, $f=7$, $g=10$, $u=3$, $v=2$, $w=1$, the  points associated with the survival of only one population, i.e. $E_1 = (3, 0, 0)$, $E_2 = (0, 2, 0)$ and $E_3 = (0, 0, 1)$, are stable, the origin $E_0 = (0, 0, 0)$ is an unstable equilibrium and the coexistence equilibrium $E_7 \approx (0.1899,    0.0270,    0.2005)$ is a saddle point.
The remaining equilibria $E_4 \approx (0.6163,0.1591,0)$, $E_5 \approx (0.2195,0,0.5317)$ and $E_6 \approx (0,0.1714,0.2647)$ are other saddle points. The manifolds joining these saddles
partition the phase space into the different domains of attraction, but intersect only at the coexistence saddle point $E_7$.

\section{Implicit partition of unity method via stable bases} \label{wsvdpu}
In this section we present the method used to reconstruct the basins of attraction. Since they are often described by implicit surfaces, 
we consider the implicit partition of unity method using locally stable bases. Such surfaces are defined by point cloud data sets of the form $ {\cal X}_N= \{ \boldsymbol{x}_i \in \mathbb{R}^{3},$ $ i=1, \ldots, N \}$, belonging to a surface in $ \mathbb{R}^{3}$. 

\subsection{Stable bases} \label{stbases}
Our goal is to recover a function $f: \Omega\to \RR$ on a bounded domain $\Omega \subset \RR^3$, sampling $f$ on $N$ pairwise 
distinct points ${\cal X}_N\subset\Omega$, namely $\boldsymbol{f} =[f_1, \ldots , f_N]^T$, $f_i = f(x_i)$, $x_i\in {\cal X}_N$. 
To do that, we consider a positive definite and symmetric kernel $\Phi : \Omega \times \Omega \rightarrow \RR$ and construct the interpolant as 
\begin{equation}\label{ansatz}
R(\boldsymbol{x}) = \sum_{j = 1}^N c_j \Phi(\boldsymbol{x}, \boldsymbol{x}_j), \qquad \boldsymbol{x}\in\Omega,
\end{equation}
where $\Phi$ is a radial kernel depending on a \textit{shape parameter} $\epsilon > 0$, 
i.e. $\Phi(\boldsymbol{x},\boldsymbol{y})=\phi_{\epsilon}(||\boldsymbol{x}-\boldsymbol{y}||_2)=\phi(\epsilon ||\boldsymbol{x}-\boldsymbol{y}||_2)$ and $\phi: \RR_{\geq 0}\to \RR$
for all $\boldsymbol{x},\boldsymbol{y} \in \Omega$. 
The coefficients $\boldsymbol{c}= [c_1, \ldots, c_N]^T$ in \eqref{ansatz} are determined 
by solving the linear system $A \boldsymbol{c}= \boldsymbol{f}$, where the interpolation 
matrix $A= [\Phi (\boldsymbol{x}_i , \boldsymbol{x}_j)]_{i,j=1}^N$. 
So the obtained solution $R$ is a function of the \textit{native Hilbert space} ${\cal N}_{\Phi}(\Omega)$ 
uniquely associated with the kernel, and, if $f\in {\cal N}_{\Phi}(\Omega)$, it is in particular 
the ${\cal N}_{\Phi}(\Omega)$-projection of $f$ into the subspace ${\cal N}_{\Phi}({\cal X}_N)$ 
spanned by the standard basis of translates ${\cal T}_{{\cal X}_N}= \{\Phi(\boldsymbol{x},\boldsymbol{x}_j), 1\leq j\leq N\}$ (cf. e.g. \cite{Fasshauer07}).

Since in most cases the matrix $A$ can be severely ill-conditioned and therefore the 
interpolant \eqref{ansatz} unstable, many efforts have recently been made to construct 
stable bases (see e.g. \cite{Fasshauer12,Fornberg07,Pazouki11}). 
In particular, here we focus on a change of basis, known as the \textit{WSVD basis}, 
described in \cite{Demarchi13, Demarchi15}. To construct this basis 
$\mathcal U = \{u_j\}_j$ of ${\cal N}_{\Phi}({\cal X}_N)$, we can assign 
an invertible coefficient matrix $D_{\mathcal U} = [d_{ij}]_{i,j=1}^N$ such that
\begin{align*}
u_j = \sum_{i=1}^N d_{ij} \Phi(\cdot, \boldsymbol{x}_i), 
\end{align*}
or, equivalently, an invertible value matrix $V_{\mathcal U} = [u_{j}(x_i)]_{i,j=1}^N$ so that 
$A = V_{\mathcal U} \cdot D_{\mathcal U}^{-1}$ \cite{Pazouki11}. 
To define the WSVD basis ${\cal U}$ for ${\cal N}_\Phi( {\cal X}_N)$, the matrices are given by 
\begin{align*}
D_{\cal U} = \sqrt{W} \cdot Q \cdot\Sigma^{-1/2}\; \mbox{ and }\; V_{\cal U} = \sqrt{W^{-1}} \cdot Q \cdot \Sigma^{1/2},
\end{align*}	
where $\sqrt{W} \cdot A \cdot \sqrt{W} =  Q \cdot \Sigma \cdot Q^{T}$ is a SVD of the weighted kernel matrix 
$A_{W}=\sqrt{W} \cdot A \cdot \sqrt{W}$, $W_{ij} = \delta_{ij} w_i$ is a diagonal matrix of positive weights, 
and $\Sigma = {\rm diag}(\sigma_1,\ldots,\sigma_N)$ the corresponding singular values. 
This basis has been constructed to mimic in a discrete sense the \textsl{eigenbasis} defined through Mercer's Theorem (see e.g. \cite{Fasshauer07}), where the inner product of $L_2(\Omega)$ is replaced with its discrete version $\ell_2({\cal X}_N)$. Moreover, as shown in \cite{Demarchi13}, the WSVD basis $\{u_j\}_{j=1}^N$ enjoys the following properties, i.e.,
\begin{enumerate}
	\item[i)]  it is ${\cal N}_{\Phi}(\Omega)$-orthonormal,
	\item[ii)]  it is $\ell_2({\cal X}_N)$-orthogonal with norm $\sigma_k$,
	\item[iii)] \label{double} $(u_k, f)_{\ell_2({\cal X}_N)} = \sigma_k(u_k, f)_{{\cal N}_{\Phi}(\Omega)}$, $\forall f \in {\cal N}_{\Phi}(\Omega)$,
	\item[iv)] $\sigma_N\geq\dots\geq \sigma_1>0$,
	\item[v)] $\sum_k\sigma_k = \phi(0)\ \mathrm{meas}(\Omega)$.
\end{enumerate}
Since the interpolation is a ${\cal N}_{\Phi}(\Omega)$-projection and thanks to property iii), we can rewrite the interpolant 
\eqref{ansatz} in terms of the ${\cal N}_{\Phi}(\Omega)$-orthonormal WSVD basis as
\begin{align*} 
R= \sum_{k=1}^{N} (f, u_k)_{{\cal N}_{\Phi}(\Omega)} u_k =\sum_{k=1}^{N} \sigma^{-1}_k (f, u_k)_{\ell_2({\cal X}_N)} u_k.
\end{align*}

\subsection{Implicit surface reconstruction via partition of unity interpolation}\label{PUMi}
To find the implicit interpolant, we need to use additional interpolation conditions considering an extra set of off-surface points. 
We construct the extra off-surface points by taking a small step away along the surface normals $\boldsymbol{n}_i$, thus obtaining 
for each data point $\boldsymbol{x}_i$ two additional off-surface points. One point lies outside the surface and is denoted by 
$\boldsymbol{x}_{N+i}=\boldsymbol{x}_i+ \delta \boldsymbol{n}_i$, whereas the other point lies inside and is denoted by 
$\boldsymbol{x}_{2N+i}= \boldsymbol{x}_i- \delta \boldsymbol{n}_i$, where $\delta$ is the stepsize \cite[Ch. 30]{Fasshauer07}. 
Now, after creating the data set, we can compute the partition of unity interpolant, whose zero contour or iso-surface interpolates 
the point cloud data ${\cal X}_N$, as well as the sets  $\cal{X}_{ \delta}^{+}=$ $  \{ \boldsymbol{x}_{N+1},\ldots, \boldsymbol{x}_{2N} \}$ and $\cal{X}_{ \delta}^{-}=$ $ \{ \boldsymbol{x}_{2N+1},\ldots,$ $\boldsymbol{x}_{3N} \}$. Some techniques to estimate the normals can be found in \cite{Hoppe94}. 

The idea of the partition of unity method is to decompose a large problem or domain $\Omega \subseteq \RR^3$ into $d$
small problems or subdomains $\Omega_j$ such that $\Omega \subseteq \bigcup_{j=1}^{d} \Omega_j$ with some mild overlap among the subdomains. 
Associated with these subdomains we construct a partition of unity, i.e. a family of compactly supported, non-negative, continuous 
functions $W_j$ with $\text{supp}(W_j) \subseteq \Omega_j$ such that 
$\sum_{j=1}^{d} W_j(\boldsymbol{x}) = 1$, $\boldsymbol{x} \in \Omega$. The global approximant thus assumes the following form
\begin{align}
\label{pui}
{\cal I}(\boldsymbol{x})= \sum_{j=1}^{d} R_j(\boldsymbol{x}) W_j(\boldsymbol{x}), \hspace{1cm} \boldsymbol{x} \in \Omega.
\end{align}
For each subdomain $\Omega_j$ we may define a {\it Shepard weight} function $W_j:\Omega_j \rightarrow \RR$ as 
\begin{align*}
W_j = \frac{\varphi_j}{\sum_{k=1}^d \varphi_k}, 
\end{align*}
$\varphi_j$ being the compactly supported Wendland $C^{2k}, \; k\ge 1$ functions \cite{Wendland05}, and a local WSVD interpolant $R_j:\Omega_j \rightarrow \RR$ of the form 
\begin{align*}
R_j = \sum_{k=1}^{N_j} \sigma^{-1}_{jk} (f, u_k^{(j)})_{\ell_2({\cal X}_j)} u_k^{(j)},
\end{align*}
where $N_j$ indicates the number of data points in $\Omega_j$, i.e. the points $\boldsymbol{x}_i^{(j)} \in {\cal X}_j = {\cal X}_N \cap \Omega_j$. 

The partition of unity approach is therefore a simple and effective computational technique because it allows us to decompose a large problem into many small subproblems, ensuring that the accuracy obtained for the local fits is carried over to the global one (for further details see \cite{Cavoretto14, Cavoretto15a, Wendland05}).

\section{Determining the basins of attraction via stable bases} \label{applicazione}
In this section we show how the method previously introduced can be used to approximate the basins of attraction of systems presenting three stable equilibria.

\begin{figure}[ht!]
	
	\includegraphics[height=.25\textheight]{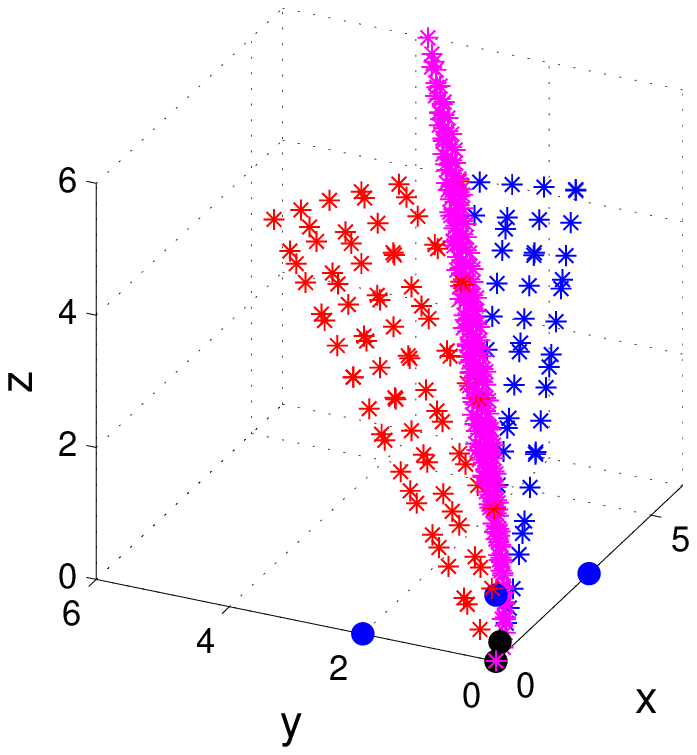} 
	\includegraphics[height=.25\textheight]{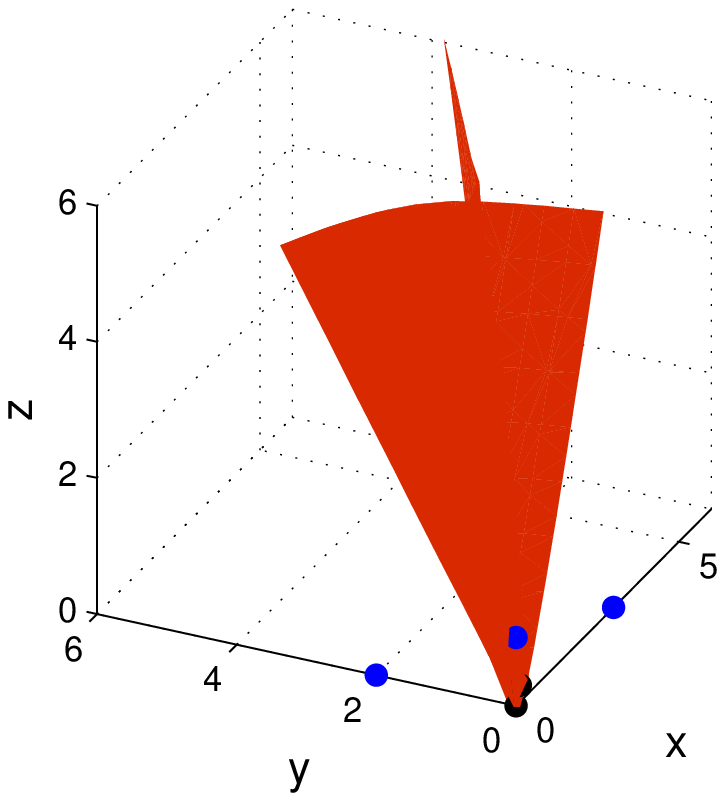} 
	\caption{Set of points lying on the surfaces determining the domains of attraction (left) and the reconstruction of the basins of attraction of $E_1$, $E_2$ and $E_3$ (right) with parameters $p = 1$, $q = 2$, $r =2$, $a = 5$, $b=4$, $c =3$, $e=7$, $f=7$, $g=10$, $u=3$, $v=2$, $w=1$. The stable equilibria are marked with the blue dot, while the unstable saddle points are represented with the black dot.}
	\label{fig1}
	
\end{figure}

In this situation we can use the routine outlined in \cite{Cavoretto15b} to approximate the basins of attraction of the system \eqref{model3d}. More precisely, we start considering $n$ equispaced points on each edge of the cube $[0,\gamma]^3$, where $\gamma \in \mathbb{R}^+$, and we define a set of initial conditions as follows
\begin{align*}
\begin{array}{llll}
P_{i_1,i_2}^{1}=(x_{i_1},y_{i_2},0) & \quad \textrm{and} & \quad P_{i_1,i_2}^{2}=(x_{i_1},y_{i_2},\gamma), & \quad i_1, i_2=1, \ldots, n, \\
P_{i_1,i_2}^{3}=(x_{i_1},0,z_{i_2}) & \quad \textrm{and} & \quad P_{i_1,i_2}^{4}=(x_{i_1},\gamma,z_{i_2}), & \quad i_1, i_2=1, \ldots, n, \\
P_{i_1,i_2}^{5}=(0,y_{i_1},z_{i_2}) & \quad \textrm{and} & \quad P_{i_1,i_2}^{6}=(\gamma,y_{i_1},z_{i_2}), & \quad i_1, i_2=1, \ldots, n.
\end{array}
\end{align*}
Then, by taking points in pairs and by performing a bisection algorithm, a certain number of separatrix points lying on the basins of attraction is found. As an example, the points lying on the separatrix manifolds shown in Figure \ref{fig1} have been found considering $n=15$ and $\gamma=6$.

In Figure \ref{fig1} we report the plot of the three surfaces describing the domains of attraction. They have been reconstructed with the stable method \eqref{pui} presented in Section \ref{wsvdpu}. Such approximate surfaces have been obtained by taking a number $d = 4$ of partition of unity subdomains and the Wendland $C^6$ function
\begin{align*}
\phi_{\epsilon}(r) = \left(1-\epsilon r\right)_+^8(32\epsilon^3r^3+25\epsilon^2r^2+8\epsilon r+1),
\end{align*} 
with shape parameter $\epsilon = 0.001$.

\end{document}